\documentstyle[11pt]{article}

\input amssym.def
\input amssym

\begin{document}

\def \a{{{\frak a}}}
\def \al{\alpha}
\def \ar{{\alpha_r}}
\def \A{{\Bbb A}}
\def \Ad{{\rm Ad}}
\def \b{{{\frak b}}}
\def \bs{\backslash}
\def \B{{\cal B}}
\def \cent{{\rm cent}}
\def \C{{\Bbb C}}
\def \CA{{\cal A}}
\def \CB{{\cal B}}
\def \CE{{\cal E}}
\def \CF{{\cal F}}
\def \CG{{\cal G}}
\def \CH{{\cal H}}
\def \CL{{\cal L}}
\def \CM{{\cal M}}
\def \CN{{\cal N}}
\def \CP{{\cal P}}
\def \CQ{{\cal Q}}
\def \CO{{\cal O}}
\def \det{{\rm det}}
\def \e{\epsilon}
\def \End{{\rm End}}
\def \Fx{{\frak x}}
\def \FX{{\frak X}}
\def \g{{{\frak g}}}
\def \ga{\gamma}
\def \Ga{\Gamma}
\def \h{{{\frak h}}}
\def \Hom{{\rm Hom}}
\def \Im{{\rm Im}}
\def \Ind{{\rm Ind}}
\def \k{{{\frak k}}}
\def \K{{\cal K}}
\def \la{\lambda}
\def \lap{\triangle}
\def \Lie{{\rm Lie}}
\def \m{{{\frak m}}}
\def \mod{{\rm mod}}
\def \n{{{\frak n}}}
\def \name{\bf}
\def \N{\Bbb N}
\def \ord{{\rm ord}}
\def \O{{\cal O}}
\def \p{{{\frak p}}}
\def \ph{\varphi}
\def \prf{{\bf Proof: }}
\def \qed{$ $\newline $\frac{}{}$\hfill {\rm Q.E.D.}\vspace{15pt}}
\def \Q{\Bbb Q}
\def \res{{\rm res}}
\def \R{{\Bbb R}}
\def \Re{{\rm Re \hspace{1pt}}}
\def \ra{\rightarrow}
\def \rank{{\rm rank}}
\def \Rep{{\rm Rep}}
\def \sign{{\rm sign}\hspace{2pt}}
\def \supp{{\rm supp}}
\def \t{{{\frak t}}}
\def \T{{\Bbb T}}
\def \tr{{\hspace{1pt}\rm tr\hspace{1pt}}}
\def \vol{{\rm vol}}
\def \V{{\cal V}}
\def \z{\frak z}
\def \Z{\Bbb Z}
\def \={\ =\ }

\newcommand{\rez}[1]{\frac{1}{#1}}
\newcommand{\der}[1]{\frac{\partial}{\partial #1}}
\newcommand{\binom}[2]{\left( \begin{array}{c}#1\\#2\end{array}\right)}
\newcommand{\norm}[1]{\parallel #1 \parallel}

\newcounter{lemma}
\newcounter{corollary}
\newcounter{proposition}
\newcounter{theorem}
\newcounter{zwisch}

\renewcommand{\subsection}{\refstepcounter{subsection}\stepcounter{lemma} 
	\stepcounter{corollary} \stepcounter{proposition}
	\stepcounter{conjecture}\stepcounter{theorem}
	$ $ \newline
	{\bf \arabic{section}.\arabic{subsection}\hspace{8pt}}}

\newtheorem{conjecture}{\stepcounter{lemma} \stepcounter{corollary} 	
	\stepcounter{proposition}\stepcounter{theorem}
	\stepcounter{subsection}\hskip-12pt Conjecture}[section]
\newtheorem{lemma}{\stepcounter{conjecture}\stepcounter{corollary}	
	\stepcounter{proposition}\stepcounter{theorem}
	\stepcounter{subsection}\hskip-7pt Lemma}[section]
\newtheorem{corollary}{\stepcounter{conjecture}\stepcounter{lemma}
	\stepcounter{proposition}\stepcounter{theorem}
	\stepcounter{subsection}\hskip-7pt Corollary}[section]
\newtheorem{proposition}{\stepcounter{conjecture}\stepcounter{lemma}
	\stepcounter{corollary}\stepcounter{theorem}
	\stepcounter{subsection}\hskip-7pt Proposition}[section]
\newtheorem{theorem}{\stepcounter{conjecture} \stepcounter{lemma}
	\stepcounter{corollary}\stepcounter{proposition}		
	\stepcounter{subsection}\hskip-11pt Theorem}[section]

\title{A simple trace formula for arithmetic groups}
\author{{\small by}\\ {} \\ Antonius H. J. Deitmar}
\date{}
\maketitle

\pagestyle{myheadings}
\markright{A SIMPLE TRACE FORMULA FOR ARITHMETIC GROUPS}

\tableofcontents

$$ $$

\begin{center} {\bf Introduction} \end{center}

Let $G$ denote a semisimple Lie group and $\Ga$ an arithmetic subgroup.
Suppose first that $\Ga$ is cocompact, i.e. the quotient manifold $G/\Ga$ is
compact.
Then the Selberg trace formula is the equality
$$
J_{geom}(f)\= J_{spectral}(f),
$$
where $J_{geom}(f)$, the  geometric side, is a sum of Orbital integrals
and the spectral side is the trace of $f$ on $L^2(G/\Ga)$.

In the case that $\Ga$ is not cocompact in $G$, the trace formula as 
developed by J. Arthur, gives a similar identity, but now both sides are 
more complicated.
The geometric side is an alternating sum of contributions belonging 
to the classes of parabolic subgroups.
This formula is set up and refined by Arthur to fit the needs of number 
theorists who primarily want to compare the trace formulas for two different 
groups.
For geometric applications one would want a trace formula whose geometric
side consists of orbital integrals just as in the cocompact case.
In this paper we show that such a formula can be deducted from
Arthur's work by plugging in special test functions.
For groups of rational rank equal to one we show how the geometric side
also can be computed further.

$$ $$

\section{A simple trace formula}\label{A simple trace formula}

We derive a simple version of Arthur's trace formula by 
inserting functions with certain restrictive properties (\ref{funct}) which 
guarantee the vanishing of the parabolic terms in the trace formula.

\subsection
Let $\CG$ be a semi simple simply connected linear algebraic group over $\Q$ and write $G$ for the group of real points.
Let $\A$ denote the adele-ring over $\Q$.
On the locally compact group $\CG(\A )$ we will fix a Haar measure.
Let $\A_{fin}$ denote the subring of finite adeles then $\A \cong\A_{fin}\times 
\R$, hence $\CG(\A ) \cong \CG(\A_{fin})\times \CG(\R)$.
We will distribute the Haar measure onto the factors.

According to Harish-Chandra \cite{HC-HA1}, every Haar-measure on $G$ comes 
from a scalar multiple $B$ of the Killing form and this form 
fixes Haar-measures on the closed subgroups on $G$.
We will keep this form $B$ at our disposal and will only fix it in later 
sections.

\subsection
We will consider $\CG$ as a subgroup of some $GL_m$.
Let $\Ga$ be a {\bf congruence subgroup}, i.e. $\Ga$ is a rational lattice which 
contains a principal congruence subgroup $\Ga (N) = \{ g\in (\CG (\Q) \cap 
GL_m(\Z)) | g\equiv 1\ \mod (N)\}$ for some natural number $N$. 

Invariantly stated the congruence property is pronounced as follows:
$\Ga \subset \CG (\Q)$ is a congruence subgroup if there is a compact open 
subgroup  $K_\Ga \subset \CG (\A_{fin})$, where $\A_{fin}$ is the ring of finite 
adeles over $\Q$, such that $\Ga = \CG (\Q) \cap K_\Ga$.
By strong approximation \cite{Kn} there is a canonical bijection
$$
\Psi :\Ga \bs G = \Ga \bs \CG(\R)\ \tilde{\longrightarrow}\ \CG(\Q)\bs 
\CG(\A )/{K_\Ga}
$$
 given by
$\Ga g \mapsto \CG(\Q) gK_\Ga$.
We will further assume $\Ga$ to be {\bf weakly neat},
 i.e. $\Ga$ is torsion free and for any 
$\ga\in\Ga$ the adjoint morphism $\Ad(\ga)$ on the $\g$ has no root of unity 
other then $1$ as an eigenvalue.
Note that any arithmetic $\Ga$ has a weakly neat subgroup of finite index 
\cite{borel}.

\subsection
Let $\CP\neq \CG$ denote a parabolic subgroup defined over $\Q$ and let 
$\CP_1\subset\CP$ be a minimal parabolic defined over $\R$. Let $\CP = \CL \CN$ 
and $\CP_1 = \CL_1 \CN_1$ denote Levi decompositions and write $\CA$ resp. 
$\CA_1$ for the split components where we assume $\CA \subset \CA_1$.
Note that we have $P=\CP(\R)=\CL(\R)^1\CA(\R)^0\CN(\R)=MAN$, where $\CL(\R)^1$ 
is the subgroup of all $m\in\CL(\R)$ such that $\chi(m)$ has absolute value $1$ 
for all rational characters $\chi$ of $\CL$.

\subsection
Fix a maximal compact subgroup $K \subset G = \CG(\R)$ such that the Lie algebra 
of $K$ is orthogonal to the Lie algebra of $\CA_1(\R)$.
Choose a maximal compact open subgroup $K_{max} = \prod_p K_p$  of 
$\CG(\A_{fin})$ such that $\CG(\A_{fin})=K_{max} \CP(\A_{fin})$, this is 
achieved by assuming $K_p$ to be a good maximal compact subgroup for all $p$.
We assume $K_\Ga \subset K_{max}$.

\subsection \label{funct}
Consider the function $f=f_{fin} \otimes f_\infty \in C_c(\CG(\A ))$ defined by
$f_{fin} = \rez{\vol(K_\Ga)}{\bf 1}_{K_\Ga}$ 
the characteristic function of the compact open subgroup divided by the volume.
On the function $f_\infty : G\ra \C$ we put the following restriction:
At first we insist that $f_\infty$ has compact support and is $j$-times 
continuously differentiable for some $j\in\N$ which is assumed to be large 
enough \cite{Art-1}.
We further insist that for any $\Q$-parabolic $P=MAN\ne G$ we have
$$
f_\infty(x^{-1}qx)=0
$$
for any $x\in G$ and $q\in MN$.

This implies that $f(man)=0$ for $man\in MAN$ if $a=1$. One might call these 
elements {\bf $P$-singular}.
Elements which are not $P$-singular will be called {\bf $P$-regular}. 
So we insist that $f_\infty$ vanishes on all $P$-singular elements for all 
nontrivial $\Q$-parabolic subgroups.
We say for short that $f$ {\bf vanishes on $\Q$-parabolically singular 
elements.}

\subsection
An element of $G$ is called {\bf regular} if its centralizer is a torus.
Regular elements are semi simple and the set $G'$ of regular elements is dense 
in 
$G$.
An element which is not regular will be called {\bf singular}.
It is not hard to show that regular elements are $P$-regular for any nontrivial 
parabolic $P$ in $G$.
This implies that if $f_\infty$ vanishes on singular elements then it already 
vanishes on parabolically singular elements.

\subsection
Given $f$ we define $K(x,y)$ by
$$
K(x,y) \= \sum_{\ga \in \CG (\Q)} f(x^{-1}\ga y).
$$
Since $f$ has compact support the sum is locally finite and so $K(x,y)$ inherits 
the 
smoothness from $f$.

\subsection
Fix a $\Q$-parabolic $\CP$.
As in \cite{Art-1}, p.923 define
$$
K_{\CP ,o}(x,y) \= \sum_{\CL (\Q)\cap o} \int_{\CN (\A )}f(x^{-1}\ga ny)\ dn.
$$
To define the geometric side, Arthur fixes some functions 
$\hat{\tau}_p(H(\delta x)-T)$.
Here the only thing we need to know is that this factor equals $1$ for the 
trivial parabolic $P=G$.
Then Arthur defines
$$
k_o^T(x,f) \= \sum_P (-1)^{\dim A_p}\sum_{\delta\in \CP(\Q)\bs\CG(\Q)}
K_{P,o}(\delta x,\delta x)\hat{\tau}_p(H(\delta x)-T),
$$
and
$$
J_{geom}(f)\= \int_{\CG(\Q)\bs \CG(\A)} \sum_\o k_o^T(x,f) dx.
$$

We will show

\begin{lemma}
Suppose $\CP\ne\CG$ and $f$ satisfies \ref{funct}, then $K_{\CP ,o}(x,x) =0$ for any 
$x\in\CG(\A )$.
\end{lemma}

\prf
Let $\ga \in \CL(\Q) \subset \CL(\A )^1$ and $n\in \CN(\A )$. Then $\ga n \in 
\CP 
(\A )^1$. 
Now let $q\in \CP (\A )^1$ be arbitrary and $x \in \CG(\A ).$ we will show that 
$f(x^{-1}qx)=0$. Assume therefore $q=q_{fin}q_\infty$ with $x^{-1}q_{fin}x \in 
\supp(f_{fin}) =K_\Ga$, then it follows that $q_{fin} \in xK_\Ga x^{-1} \cap 
\CP(\A )$, a compact subgroup of $\CP (\A )$. 
Any continuous quasicharacter with values in $]0,\infty [$ will therefore be 
trivial on $q_{fin}$, hence $q_{fin} \in \CP (\A )^1$. Since $q$ already was in 
$\CP(\A )^1$ it follows $q_\infty \in \CP(\A )^1 \cap \CP(\R) = 
\CL_Q'(\R)\CN_Q(\R)=MN$ but this implies $f_\infty(q_\infty)=0$ as claimed.
\qed

\subsection
For $g\in C_c(G)$ and $y\in G$ define the {\bf orbital integral} as
$$
\CO_y(g)\ :=\ \int_{G_y\bs G}g(x^{-1} yx)\ dx.
$$ 
Where $G_y$ denotes the centralizer of $y$ in $G$.
(Recall our conventions on Haar-measures.)
It is known that the integral always converges.

\subsection
Let $R$ denote the representation of $\CG(\A )$ on $L^2(\CG(\Q)\bs \CG(\A ))$ 
and 
let $R(f)$ denote the operator defined by $f$. 
From the lemma it follows that along the diagonal the kernel $K$ coincides with 
the modified kernel as in \cite{Art-1}.
So it follows that the integral over $\CG(\Q)\bs \CG(\A )$ of the diagonal 
$K(x,x)$ exists.

\begin{theorem}
Let the function $f$ on $G(\A)$ satisfy \ref{funct} then the geometric 
side of the trace formula $J_{geom}(f)$ equals:
$$
\int_{\CG (\Q) \bs \CG(\A )} K(x,x)\ dx \= \sum_{[\ga]}\ \vol (\Ga_\ga \bs 
\CG_\ga(\R))\ \CO_\ga(f_\infty),
$$
where the sum on the right hand side runs over the set of all conjugacy classes 
$[\ga]$ in the group $\Ga$.
\end{theorem}

\prf
Consider the bijection $\Psi : \Ga \bs G \ra \CG(\Q) \bs \CG(\A ) /K_\Ga$ and 
let 
$\Psi_*$ denote the unitary map
$$
\Psi_* : L^2(\Ga \bs G) \ra L^2(\CG(\Q)\bs \CG(\A ) /K_\Ga ) 
\begin{array}{c}\sqrt{\vol(K_\Ga)}^{-1}\\ {\hookrightarrow}\\ {}\end{array} 
L^2(\CG(\Q)\bs \CG(\A )),
$$
further let
$$
\Phi_* : L^1(\Ga \bs G) \ra L^1(\CG(\Q)\bs \CG(\A ) /K_\Ga ) 
\begin{array}{c}{\vol(K_\Ga)}^{-1}\\ {\hookrightarrow}\\ {}\end{array}
L^1(\CG(\Q)\bs \CG(\A )).
$$
At first note that by construction
$$
\int_{\CG(\Q)\bs \CG(\A )} \Phi_*\ph(x) dx \= \int_{\Ga\bs G}\ph(x) dx.
$$
Next we have $\Psi_* R(f_\infty)=R(f)\Psi_*$ and this gives 
$\Psi_{*,1}\Psi_{*,2}K =K_\infty$, where $R(f_\infty)$ is the convolution by 
$f_\infty$ on $L^2(\Ga \bs G)$ and $K_\infty$ its kernel.
The indices at $\Psi_*$ indicate that it is applied to each argument separately.
Let $K(x) =K(x,x)$ and $K_\infty(x)=K_\infty(x,x)$ as a function in one 
argument.
Clearly $(\Psi_{*,1}\Psi_{*,2}K)(x,x) = \Phi_*K(x)$ and thus
$$
\int_{\CG(\Q) \bs \CG(\A )} K(x,x)\ dx \= \int_{\Ga \bs G}K_\infty (x,x)\ dx.
$$
Arthur showed absolute convergence in \cite{Art-1}.
The rest is the usual calculation expressing the integral of the diagonal as a 
sum of orbital integrals.
\qed

\subsection
Now we consider the spectral side. 
We will now assume that the rank of $\CG$ over $\Q$ equals $1$.
Then there is, up to conjugation, only one nontrivial $\Q$-parabolic 
$\CP$ in $\CG$.
By Arthur's kernel identity we get
that 
$$
\int_{\CG (\Q) \bs \CG(\A )} K(x,x) dx
$$
equals the sum of 
$$
\sum_{\pi \in \hat{\CG (\A )}}m_\pi \tr\ \pi(f)
$$
and
$$
\sum_\chi \rez{4\pi i} \int_{i\a_{_0}^*}\sum_{\phi \in \CB_{\CP ,\chi}} 
\left( \wedge^T E(. ,I_\CP (\la ,f)\phi,\la),\wedge^TE(.,\phi,\la)
\right) d\la,
$$
where $m_\pi$ is the multiplicity of $\pi$ in 
$L^{2}(\CG (\Q)\bs \CG(\A ))_{disc}$, the discrete part, further 
the sum $\sum_\chi$ runs over all cuspidal representations of 
$\CL(\A )$ and $\CB_{\CP,\chi}$ is a orthonormal basis of the 
$\chi$-isotype and $\wedge^T E$ is the truncated Eisenstein series 
as in \cite{Art-Corv}. 
Using the fact that the above is independent of the cutoff parameter $T$, 
by the Maa\ss -Selberg relations (\cite{Langl},\cite{Art-Corv}) and 
the usual Fourier analysis \cite{GelbJac} the latter summand is seen to be
$$
 \rez{4\pi i} \int_{i\a_{_0}^*} \tr (1+M(\la)^{-1}M'(\la)  
I_\la(f)) d\la\ -\ \rez{4} \tr (M(0)I_0(f)),
$$
where $M(\la)$ is the intertwining operator of \cite{Art-Corv} and 
$I_\la$ is the induced representation on the space of which $M(\la )$ acts.
For the convenience of the reader we will recall the definition of $I_\la$ and 
$M(\la)$ next.

\subsection
We have the decomposition $\CG(\A ) = \CN(\A ) \CL(\A )^1\CA(\A )^0K_{max}K$.
In this decomposition any $y\in\CG(\A )$ writes as $y=nm\exp(Y) k_fk$ where 
$Y\in\a$ is uniquely determined.
We define a map $H_P : \CG(\A )\ra \a$ mapping $y$ to $Y$.

Let $\CH_P^0$ be the space of functions
$$
\ph : _{\CN(\A )\CL(\Q)\CA(\R)^0}\bs^{\CG(\A )}\ \ra\ \C
$$
such that
\begin{itemize}
\item
for any $x\in \CG(\A )$ the function $m\mapsto \ph(mx)$, $m\in\CL(\A )$ is 
$\z_{\CL(\R)}$-finite, where $\z_{\CL(\R)}$ is the center of the universal 
enveloping algebra of $\l ={\rm Lie}(\CL(\R))$,
\item
$\ph$ is right $K_{max}K$-finite,
\item
$$
\norm{\ph}^2 = \int_{K_{max}K}\int_{\CA(\R)^0\CL(\Q)\bs\CL(\A )}\mid 
\ph(mk)\mid^2 dmdk <\infty .
$$
\end{itemize}
The last point defines an obvious scalar product on $\CH_P^0$.
Let $\CH_P$ be the Hilbert space completion of $\CH_P^0$.
For $\la\in\a$, $\ph\in\CH_P$ and $x,y\in\CG(\A )$ put
$$
(I_\la(y)\ph)(x) \= \ph(xy)e^{\langle \la +\rho ,H_P(xy)-H_P(x)\rangle}.
$$
We will write $I_\la$ for the resulting $\CG(\A )$-representation.
Note that this is just the induced representation from $\CP(\A )$ to $\CG(\A )$ 
of 
$L^2(\CL(\Q)\bs \CL(\A )^1)\otimes\la\otimes 1$.

\subsection
Let $w\in K_{max}K \cap\CG(\Q)$ be a representative for the nontrivial element 
of the Weyl group $W(\g ,\a)$.
For $\la\in\a$ with $\Re(\la)>>0$ we define an operator $M(\la)$ on $\CH_P^0$ by
$$
M(\la)\ph(x) = \int_{\CN(\A )} \ph(wnx) e^{\langle \la +\ph ,H_P(wnx)\rangle 
+\langle \la -\rho ,H_P(x)\rangle} dn.
$$
Langlands has shown that this operator extends meromorphically in $\la$ and 
satisfies
$$
M(-\la) =M(\la)^{-1},\ \ \ M(\la)^* = M(\bar{\la}).
$$
These equations especially imply that $M(\la)$ is unitary if $\la$ is purely 
imaginary.

Let $\Pr : \CH_P\ra\CH_P^{K_\Ga}$ be the projection to the space of 
$K_\Ga$-fixed vectors.
In our considerations  only the operator $M(\la)^\Ga :=\Pr M(\la)\Pr$ will 
occur.

\subsection \label{induced_rep}
Since $\CG$ has $\Q$-rank one it follows that $\CL(\Q)\bs \CL(\A )^1$ is 
compact.
So $L^2(\CL(\Q)\bs \CL(\A )^1)$ decomposes as $\CL(\A )^1$-representation
$$
L^2(\CL(\Q)\bs \CL(\A )^1) \=\bigoplus_{\chi} N(\chi)\chi ,
$$
where the sum runs over the unitary dual of $\CL(\A )^1$ and the multiplicities 
$N(\chi)$ are finite.
Therefore
\begin{eqnarray*}
I_\la & = & {\rm Ind}_{\CP(\A )}^{\CG(\A )}(L^2(\CL(\Q)\bs \CL(\A )^1) 
\otimes\la\otimes 1)\\
	&=& \bigoplus_\chi N(\chi) {\rm Ind}_{\CP(\A )}^{\CG(\A )}(\chi 
\otimes\la\otimes 1)
\end{eqnarray*}
Write $\pi_{\chi ,\la}$ for ${\rm Ind}_{\CP(\A )}^{\CG(\A )}(\chi 
\otimes\la\otimes 1)$.
We consider $\chi\otimes\la =\chi_\la$ as an irreducible admissible  
representation of $\CL(\A )=\CL(\A )^1\CA(\R)^0$.
Since any irreducible admissible representation of $\CL(\A )$ can 
be written as a tensor product of local representations we get
$$
\chi_\la \=\bigotimes_p \chi_{\la ,p}\otimes \chi_{\la ,\infty}.
$$
For example, $\chi_{\la ,\infty}$ is a representation of $\CL(\R) =MA$.
For later use we introduce the notation $\chi_\infty := \chi_{\la ,\infty}|_M$, 
which does not depend on $\la$.
Accordingly we get
$$
\pi_{\chi_\la} \= \bigotimes_p \pi_{\chi_\la ,p} \otimes \pi_{\chi_\la ,\infty} 
\= \pi_{\chi_\la , fin} \otimes \pi_{\chi_\la ,\infty}.
$$
Note that with this notation we have 
$\pi_{\chi_\la ,\infty}=\pi_{\chi_\infty ,\la}$.
For the space of $K_\Ga$-fixed vectors we get
$$
(\pi_{\chi_\la})^{K_\Ga} \= (\pi_{\chi_\la , fin})^{K_\Ga} \otimes \pi_{\chi_\la 
,\infty}.
$$
The space $(\pi_{\chi_\la , fin})^{K_\Ga}$ is finite dimensional, its dimension 
$N_\Ga(\chi)$ does not depend on $\la$.
As a $G$-representation
$$
I_\la^{K_\Ga} \= \bigoplus_\chi N(\chi) N_\Ga(\chi) \pi_{\chi_\la ,\infty}.
$$

The representation $\pi_{\chi_\la ,\infty}$ equals the direct Hilbert sum of its 
$K$-isotypes
$$
\pi_{\chi_\la ,\infty} \= \bigoplus_{\tau\in\hat{K}}\pi_{\chi_\la 
,\infty}(\tau),
$$
and each isotype is finite dimensional.
The Weyl-group $W=W(A,G)$ acts on the unitary dual $\hat{M}$ of the group
$M=\CL(\R)^1$.
Write $\CO$ for an orbit in $\hat{M}$, then
$$
\CH_P^{K_\Ga} \= \bigoplus_\CO \bigoplus_{\tau} \pi_\CO(\tau),
$$
with $\pi_\CO(\tau) = \bigoplus_{\chi_{\la ,\infty}\in\CO} \pi_{\chi_{\la 
,\infty}}(\tau)$ (here any $\la$ will give the same space, but the 
representation on it varies.)
The spaces $\pi_\CO(\tau)$ are finite dimensional and stable under $M(\la)$ for 
any $\la$.

\subsection \label{kern_ist_PW}
For any measure space $(\Omega ,\mu)$ and any Hilbert space 
$V$ we write $L^2(\Omega ,V)$ for the space of all $V$-valued 
square integrable functions on $\Omega$.

We now come to the discussion of 
$\pi_{\chi ,\la}(f) = Pr_\Ga \otimes \pi_{\chi_\infty ,\la}(f_\infty)$; recall 
\ref{induced_rep} for notations.
Using the compact model of induction we see that $\pi_{\chi_\infty 
,\la}(f_\infty)$ 
can be written as an integral operator on the space
$$
L^2(K ,\chi_\infty) := \{ \ph \in L^2(K ,V_{\chi_\infty})| \ph(mk)
\=\chi_\infty(m)f(k),\ m\in K \cap M\},
$$
with kernel
\begin{eqnarray*}
k_{f,\la}(k,k') &=& \int_{MAN} f_\infty (k^{-1} mank') a^{\rho +\la} 
\chi_\infty(m) dm da dn\\
	&=& \int_\R \tilde{f}(k,k',\exp (tH)) e^{\la(H) t} dt,
\end{eqnarray*}
for a compactly supported function $\tilde{f}$ on $K\times K\times A$.
Thus we may view the kernel $k_{f,\la}$ pointwise as a Paley-Wiener function in 
$\la$.

\begin{lemma}
For any irreducible admissible representation $\chi_\la$ of $\CL(\A )$ it holds
$$
\int_{i\a_0^*} \tr \pi_{\chi_\la ,\la}(f) d\la \= 0.
$$
\end{lemma}

\prf
Since $\pi_{\chi_\la}(f) = Pr_\Ga \otimes \pi_{\chi_\la ,\infty}(f_\infty)$ it 
remains to show
$$ 
\int_{i\a_0} \tr \pi_{\chi_\la ,\infty}(f_\infty) d\la \= 0.
$$
Let $V_{\chi_\infty} =\bigoplus_{\tau\in\hat{K_M}}V_{\chi_\infty}(\tau)$ be the 
decomposition of $V_{\chi_\infty}$ into $K_M=K \cap M$-isotypes.
To this corresponds the decomposition
$$
L^2(K ,\chi_\infty) \=\bigoplus_{\tau\in\hat{K_M}} L^2(K ,\chi_\infty ,\tau),
$$
where
$$
L^2(K ,\chi_\infty ,\tau) \= \{ \ph\in L^2(K ,V_{\chi_\infty}(\tau)) | 
\ph(mk)=\chi_\infty(m)f(k),\ m\in K_M\}.
$$

Let $Pr_\tau$ denote the projection $L^2(K ,\chi_\infty)\ra L^2(K ,\chi_\infty 
,\tau)$.
Since the operator $\pi_{\chi_\la ,\infty}(f)$ is given by the smooth kernel 
$k_f(k,k')$ it follows that the operator $\pi_{\chi_\la ,\infty}(f) Pr_\tau$ is 
given by the kernel
$$
k_{f,\tau}(k,k') \= \int_{K_M} k_f(k,k_Mk'){\tr\tau(k_M)} dk_M.
$$
The latter is smooth, so it follows
$$
\tr(\pi_{\chi_\la ,\infty}(f_\infty)Pr_\tau) \= \int_K \tr(k_{f,\tau}(k,k)) dk
$$ $$
\= \int_\R \int_K \tr \tilde{f}_\tau (k,k,\exp(tH))e^{t\la (H)} dt,
$$
where $\tilde{f}_\tau (k,k',a)$ is defined to be
$$
\int_{K_M} \int_{\CL(\R)^1 \CN(\R)} f_\infty (k^{-1} mank_Mk') e^{t\la(H)} 
\chi_\infty(m) {\tr\tau(k_M)}dm dn dk_M.
$$

It follows
$$
\int_{i\a_0} \tr\pi_{\chi_\la ,\infty}(f_\infty) d\la \= \sum_\tau 
\int_{i\a_0}\int_\R\int_K \tr\tilde{f}_\tau(k,k,\exp(tH)) e^{t\la(H)} d\la
$$ $$
\= \rez{2\pi} \sum_\tau \int_K \tr \tilde{f}_\tau (k,k,1) dk
$$
by the Fourier-inversion theorem.
Since $f_\infty$ vanishes on $P$-singular elements it follows that 
$\tilde{f}_\tau (k,k,1)=0$ for any $k\in K$, $\tau\in\hat{K_M}$.
\qed 

Up to this point we have shown:

\begin{theorem} \label{simple_tf_1}
Let the function $f$ on $\CG(\A )$ satisfy \ref{funct}, then the sum
$$
\sum_{[\ga]} \vol(\Ga_\ga\bs G_\ga) \CO_\ga(f_\infty)
$$
equals
$$
\sum_{\pi\in\hat{G}}N_\Ga(\pi)\tr\pi(f_\infty) \ -\ \rez{4}\tr(M(0) 
I_0(f))
$$ $$
+\  \rez{4\pi i} \int_{i\a_0^*} \tr((M(\la)^{-1}M'(\la)) 
I_\la(f)) d\la .
$$
\end{theorem}

\section{The continuous contribution}\label{The continuous contribution}
In this section we will give the continuous contribution
$$
\rez{4\pi i} \int_{i\a_0^*} \tr((M(\la)^{-1}M'(\la)) 
I_\la (f)) d\la
$$ 
a different shape.
This section is more general then the rest of the paper since we can take for 
$f_\infty$ an arbitrary element of $C_c^j(G)$.

\subsection
For $r>0$ and $a\in\C$ let $B_r(a)$ be the closed disk around $a$ of radius $r$.
Let $g$ be a meromorphic function on $\C$ with poles $ a_1 ,a_2,\dots $.
We say that $g$ is {\bf essentially of moderate growth} if there is a natural 
number $N$, a constant $C>0$ and a sequence of positive real numbers $r_n$, 
tending to zero such that the disks $B_{r_n}(a_n)$ are pairwise disjoint and 
that 
on the domain $D=\C -\bigcup_n B_{r_n}(a_n)$ it holds
$|g(z)|\ \le\ C|z|^N.$
In that case the constant $N$ is called the {\bf growth exponent}.

\begin{lemma} \label{logarithmic_derivative_essent_mod_grow}
Let $f$ be a meromorphic function of finite order and let $g=f'/f$ be its 
logarithmic derivative.
Then $g$ is essentially of moderate growth with growth exponent equals the order 
of $f$ plus two.
\end{lemma}

\prf
This is a direct consequence of Hadamard's factorization theorem applied to $f$.
\qed

\subsection
Recall the definition of $M(\la)$.
The integral over $\CN(\A )$ may be written as a product of an integral over 
$\CN(\A_{fin})$ and an integral over $\CN(\R)$.
Thus $M(\la) =M_{fin}(\la)\otimes M_\infty(\la)$, and so 
$M(\la)^{-1}M'(\la)=M_{fin}(\la)^{-1}M_{fin}'(\la) + 
M_\infty(\la)^{-1}M_\infty'(\la)$.
According to \cite{Art-intresI}, Theorem 2.1 we write 
$M_{fin}(\la)|_{\pi_\CO(\tau)} = r_{fin}(\la)R_{fin}(\la)$ with a scalar-vlued 
meromorphic function $r_{fin}(\la)$.
Property $(R_6)$ in Theorem 2.1 of loc. cit. implies that $R_{fin}(\la)$ is of 
finite order, the order being independent on $\CO$ and $\tau$, hence the same 
holds for $\det(R_{fin}(\la))$.
Further p. 39 of loc. cit. shows that $r_{fin}$ is of finite order independent 
on $\CO$ and $\tau$.

\begin{lemma} \label{ess_mod_gro1}
The function $\la \mapsto \tr M_{fin}(\la)^{-1}M_{fin}'(\la)|_{\pi_\CO(\tau)}$ 
is of essentially moderate growth with growth exponent independent of $\CO$ and 
$\tau$.
\end{lemma}

\prf
Let $\psi(\la) := \det(M_{fin}(\la)|_{\pi_\CO(\tau)})$ then $\tr 
M_{fin}(\la)^{-1}M_{fin}'(\la)|_{\pi_\CO(\tau)}$ equals $\psi'/\psi(\la)$, so by 
Lemma \ref{logarithmic_derivative_essent_mod_grow} it suffices to show that the 
order of $\psi$ is independent on $\CO$ and $\tau$.
This is clear by the above.
\qed

\begin{lemma} \label{ess_mod_gro2}
The matrix valued function $M_\infty(\la)^{-1}M_\infty'(\la)|_{\pi_\CO (\tau)}$ 
is essentially of moderate growth with growth exponent independent on $\CO$ and 
$\tau$.
\end{lemma}

\prf 
For the length of this proof write $M(\la)$ for $M(\la)|_{\pi_\CO(\tau)}$ and 
$M_\infty(\la)$ for $M_\infty(\la)|_{\pi_\CO(\tau)}$.
In \cite{mu} ,p. 514, it is shown that $M(\la)$ is a matrix valued meromorphic 
function of order $\dim G/K +2$.
Above we showed that $M_{fin}$ is of finite order independent of $\CO$ and 
$\tau$.
Together the same follows for $M_\infty(\la)$.
According to Theorem 2.1 in \cite{Art-intresI} we have 
$M_\infty(\la)=r_\infty(\la)R_\infty(\la)$ and hence 
$M_\infty(\la)^{-1}M_\infty'(\la) =r_\infty'/r_\infty(\la) + 
R_\infty(\la)^{-1}R_\infty'(\la)$.
By formula (3.5) of \cite{Art-intresI} it follows that the normalizing factor 
$r_\infty$ satisfies the same growth conditions as $M_\infty(\la)$.
Now Lemma \ref{logarithmic_derivative_essent_mod_grow} applies to the first 
summand.
The second summand is rational by p. 29 in \cite{Art-intresI} and the proof on 
page 37 of \cite{Art-intresI} implies that the degree of $R_\infty(\la)$ and 
$R_\infty(\la)^{-1}$ and hence of $R(\la)^{-1}R'(\la)$ does only depend on $G$.
\qed

\subsection
We want to give the integral over $i\a_0^*$ in Theorem \ref{simple_tf_1} a 
different 
shape. 
To this end recall that the kernel $k_{f,\la}$ is a Paley-Wiener function in the 
argument $\la$.

We will formulate a general remark on Paley-Wiener functions.
For a natural number $n$ let $C_c^n(\R)$ denote the space of $n$-times 
continuously differentiable compactly supported functions on $\R$.
By a {\bf Paley-Wiener function of order $n$} we mean a function $h$ which is 
the Fourier transform of some $g\in C_c^n(\R)$.
Since it better fits into our applications we will change coordinates from $z$ 
to $iz$. So a Paley-Wiener function $h$ will be of the form
$$
h(z) \= \int_{-\infty}^\infty g(t) e^{zt} dt
$$
for some $g\in C_c^n(\R)$.

\begin{proposition}
Let $h$ be a Paley-Wiener function of order $n$ and fix $a\in \C$.
There is a unique decomposition
$$
h \= h_a^{+,n} + h_a^{-,n}
$$
such that the functions $h^{\pm n}_a$ are holomorphic in $\C -\{ a\}$, both have 
at most a pole of order $<n$ at $a$. 
Further for some $C>0$ the following estimates hold:
\begin{eqnarray*}
|h_a^{+,n}(z)|&\leq & \frac{C}{|z-a|^n}\ \ \ {\rm for}\ \Re (z) \leq 0, \ z\neq 
a,\\
|h_a^{-,n}(z)|&\leq & \frac{C}{|z-a|^n}\ \ \ {\rm for}\ \Re (z) \geq 0, \ z\neq 
a.
\end{eqnarray*}
\end{proposition}

\prf
Let us show uniqueness first.
Suppose we are given two decompositions $h=h^++h^-=h_1^++h_1^-$ of the above 
type then $\tilde{h}=h^+-h_1^+=h_1^--h^-$ satisfies 
$|\tilde{h}(z)|\le\frac{2C}{|z-a|^n}$ for all $z\ne a$. 
Therefore the entire function $(z-a)^n\tilde{h}(z)$ is bounded, hence constant.
But this function vanishes at $a$ by the pole order condition, whence the claim.

For the existence assume
$$
h(z) \= \int_{-\infty}^\infty g(t) e^{zt} dt
$$
for some $g\in C_c^n(\R)$.
Now define
$$
h_a^{+,n}(z) \ :=\ (\frac{1}{z-a})^n \int_0^\infty 
(g(t)e^{at})^{(n)} e^{(z-a)t} dt -\frac{c(g)}{(z-a)^n}
$$
and
$$
{h}_a^{-,n}(z) \ :=\ (\frac{1}{z-a})^n \int_{-\infty}^0 
(g(t)e^{at})^{(n)} e^{(z-a)t} dt +\frac{c(g)}{(z-a)^n},
$$
where $c(g)=\int_0^\infty (g(t)e^{at})^{(n)} dt$.
Partial integration shows that ${h} = {h}_a^{+,n} +{h}_a^{-,n}$, the 
rest is clear.
\qed

\subsection \label{orth_pol}
Note that if $g$ vanishes at $t=0$ to order $j+1$ and $n\le j$, then
$$
{h}^{\pm ,n}_a \= {h}^{\pm ,n-1}_a \=\dots \= {h}^{\pm ,1}_a 
$$
and this further equals
$$
h^\pm(z) \ :=\ \int_0^\infty g(\pm t) e^{\pm tz} dt.
$$
In this case we say that $h$ is {\bf orthogonal to polynomials of degree $\le 
j$}.

If $a=0$, we will generally drop the index, so ${h}^{\pm ,n}_0={h}^{\pm ,n}$.

Finally note that, by the formula given above one sees that if $g$ depends  
differentiably or holomorphically on some parameter then the same holds for 
$h^\pm$.

\subsection \label{residuensatz_fuer_verflechtungsoperator}
Fix some $n\leq j$, but still large
 and denote by $k_{f,\la ,a}^{\pm ,n}$ the kernels we get by applying this 
construction to $k_{f,\la}$ as a function in $\la$.
Write $T_{f,\la ,a}^{\pm ,n}$ for the corresponding operator at infinity 
and $I_{\la ,a}^{\pm ,n}(f)$ for the global operator 
$Pr_\Ga \otimes T_{f,\la ,a}^{\pm ,n}$.
Suppose $a\in\a^*$ has negative real part and does not coincide with a 
pole of $M(\la)^{-1}M'(\la)$.
We get that
$$
\rez{4\pi i} \int_{i\a_0^*} \tr M(\la)^{-1}M'(\la)I_{\la}(f) d\la
$$
equals
\begin{eqnarray*}
&{}&\rez{4\pi i} \int_{i\a_0^*} \tr M(\la)^{-1}M'(\la)I_{\la 
,a}^{+,n}(f) d\la\\
&+& \rez{4\pi i} \int_{i\a_0^*} \tr M(\la)^{-1}M'(\la)I_{\la 
,a}^{-,n}(f) d\la.
\end{eqnarray*}
We move the integration paths to the left and the right resp. to get the 
residues plus a term which  tends to zero according to the Lemmas 
\ref{ess_mod_gro1} and \ref{ess_mod_gro2}.
The above becomes
\begin{eqnarray*}
& &\rez{2} \sum_{\Re \la <0} \tr R_\la I_{\la ,a}^{+,n}(f)
\\
&+& \rez{2} \res_{\la =a} \tr M(\la)^{-1}M'(\la)I_{\la ,a}^{+,n}(f)
\\
&-& \rez{2} \sum_{\Re \la >0} \tr R_\la I_{\la ,a}^{-,n}(f),
\end{eqnarray*}
where $R_{\la_0} := \res_{\la =\la_0} M(\la)^{-1}M'(\la)$.

\subsection
We say that a function $f\in C_c^j(G)$ is {\bf orthogonal to polynomials of 
degree $\le j$} if the operator valued function $\la \mapsto \pi_{\xi ,\la}(f)$ 
satisfies this condition for any $\xi\in\hat{M}$.
In that case \ref{orth_pol} immediately gives that the above equals
$$
\rez{2} \sum_{\Re \la <0} \tr R_\la I_\la^+(F) -\rez{2} \sum_{\Re \la >0} \tr 
R_\la I_\la^-(F).
$$

Furthermore $f$ is called {\bf even} if $\la \mapsto \pi_{\xi ,\la}(f)$ is an 
even function for any $\xi\in\hat{M}$.
In that case the functional equation gives $\tr R_{-\la}I_{-\la}^-(f) =\tr R_\la 
I_\la^+(f)$ and so for $f$ even and orthogonal to polynomials of degree $\le j$ 
we end up with the simple expression
$$
\sum_{\Re \la <0} \tr R_\la I_\la^+(f).
$$

\subsection
Recall the definition of $M(\la)$ as an integral over $\CN(\A )$.
This integral can be written as a product of an integral over $\CN(\A_{fin})$ 
and an integral over $\CN(\R)$ giving a decomposition 
$M(\la)=M_{fin}(\la)\otimes M_\infty(\la)$.
The second factor, restricted to the contribution of a single 
$\chi\in\hat{\CL(\A )^1}$, coincides with the Knapp-Stein intertwining operator 
at 
the infinite place and can separately be continued to a meromorphic function on 
the plane \cite{Wall2}.
Therefore the first also extends meromorphically and we get
$M(\la)^{-1}M'(\la) = M_{fin}(\la)^{-1}M_{fin}'(\la) + 
M_\infty(\la)^{-1}M_\infty'(\la)$ and so $R_\la = R_{fin ,\la}\otimes 1 + 
1\otimes R_{\infty ,\la}$ which implies that $\tr R_\la \pi_\la^{+,n}(f)$ equals
$$
\sum_\chi N(\chi)\tr (R_{fin ,\chi_\la}^\Ga) \tr\pi_{\chi_\la 
,\infty}^{+,n}(f_\infty)
+ \sum_\chi N(\chi)N_\Ga(\chi) \tr R_{\infty ,\chi_\la}\pi_{\chi_\la 
,\infty}^{+,n}(f_\infty).
$$

\begin{lemma}
The number $N_\Ga(\chi ,\la) := \tr (R_{fin ,\chi_\la}^\Ga)$ is an integer.
\end{lemma}

\prf
Let $A(\la)$ be the finite dimensional matrix 
$M_{fin ,\chi}(\la)^\Ga$, then
\begin{eqnarray*}
\tr (R_{fin ,\chi_{\la_0}}^\Ga) &=& \tr (\res_{\la =\la_0} A(\la)^{-1}A'(\la))\\
	&=& \res_{\la =\la_0} \tr A(\la)^{-1}A'(\la)\\
	&=& \res_{\la =\la_0} \tr \der{\la} \log A(\la)\\
	&=& \res_{\la =\la_0} \der{\la} \log \det A(\la)
\end{eqnarray*}
Let $c(\la)=\det(A(\la))$ then the last expression becomes
$\res_{\la =\la_0} \frac{c'(\la)}{c(\la)}$,
which is an integer.
\qed

\subsection \label{fourier_trafo}
Let $\hat{f}_\infty$ be the trace of the Fourier-transform of $f_\infty$, i.e. 
for an admissible representation $\pi$ of $G$ of finite length let
$\hat{f}_\infty(\pi) \ :=\ \tr\pi(f_\infty).$
The function $\la\mapsto\hat{f}_\infty(\pi_{\chi_\la ,\infty})$ is a 
Paley-Wiener function and so it decomposes as above where we write
$\hat{f}_\infty(\pi_{\chi_\la ,\infty})\= \hat{f}_\infty^{+,n}(\chi 
,\la)+\hat{f}_\infty^{-,n}(\chi ,\la).$
Since the trace equals the integral over the kernel it follows
$\tr\pi_{\chi_\la ,\infty}^{+,n}(f_\infty) \= \hat{f}_\infty^{+,n}(\chi ,\la).$ 
We can summarize the results of this section in

\begin{proposition} \label{cont_cont}
Let $f_\infty$ in $C_c^j(G)$ be even and orthogonal to polynomials of degree 
$\le j$, then
$$
\rez{4\pi i} \int_{i\a_0} \tr((M(\la)^{-1}M'(\la)) 
I_\la(f)) d\la
$$ 
equals
$$
\sum_{\Re \la <0} \tr R_\la I_\la^+(f),
$$
and this can be written as
$$
+ \sum_\chi N(\chi) \left( \sum_{\Re(\la)<0} N_\Ga(\chi ,\la) 
\hat{f}^{+}_\infty(\chi ,\la)+N_\Ga(\chi)\tr R_{\infty 
,\chi_\la}\pi^{+}_{\chi_\la ,\infty}(f_\infty)\right) .
$$
\end{proposition}


\small
\hspace{-20pt}
School of Mathematical Sciences\\
University of Exeter\\
Laver Building, North Park Road\\
Exeter EX4 4QE\\
Devon, UK

\end{document}